%

\magnification=\magstep1
\def\forces{\parallel\!\!\! -}


\def\hexnumber#1{\ifcase#1 0\or1\or2\or3\or4\or5\or6\or7\or8\or9\or
	A\or B\or C\or D\or E\or F\fi }

\font\teneuf=eufm10
\font\seveneuf=eufm7
\font\fiveeuf=eufm5
\newfam\euffam
\textfont\euffam=\teneuf
\scriptfont\euffam=\seveneuf
\scriptscriptfont\euffam=\fiveeuf
\def\frak{\fam\euffam \teneuf}

\font\tenmsx=msam10
\font\sevenmsx=msam7
\font\fivemsx=msam5
\font\tenmsy=msbm10
\font\sevenmsy=msbm7
\font\fivemsy=msbm5
\newfam\msxfam
\newfam\msyfam
\textfont\msxfam=\tenmsx  \scriptfont\msxfam=\sevenmsx
  \scriptscriptfont\msxfam=\fivemsx
\textfont\msyfam=\tenmsy  \scriptfont\msyfam=\sevenmsy
  \scriptscriptfont\msyfam=\fivemsy
\edef\msx{\hexnumber\msxfam}

\mathchardef\upharpoonright="0\msx16
\let\restriction=\upharpoonright
\def\Bbb#1{\tenmsy\fam\msyfam#1}

\def\re{{\restriction}}
\def\omom{{\omega^\omega}}
\def\omlom{{\omega^{< \omega}}}

\def\twolom{{2^{<\omega}}}

\def\Smallskip{\vskip1.4truecm}
\def\Bigskip{\vskip2.2truecm}
\def\Hoskip{\hskip1.6truecm}

\def\qed{{\vcenter{\hrule height.4pt \hbox{\vrule width.4pt height5pt
 \kern5pt \vrule width.4pt} \hrule height.4pt}}}
\def\ok{\vbox{\hrule height 8pt width 8pt depth -7.4pt
    \hbox{\vrule width 0.6pt height 7.4pt \kern 7.4pt \vrule width 0.6pt height 7.4pt}
    \hrule height 0.6pt width 8pt}}
\def\nt{{\leq}\kern-1.5pt \vrule height 6.5pt width.8pt depth-0.5pt \kern 1pt}
\def\sd{{\times}\kern-2pt \vrule height 5pt width.6pt depth0pt \kern1pt}
\def\zp#1{{\hochss Y}\kern-3pt$_{#1}$\kern-1pt}

\def\bb{{\frak b}}

\def\dd{{\frak d}}
\def\ee{{\frak e}}
\def\pp{{\frak p}}
\def\se{{\frak{se}}}
\def\ss{{\frak s}}

\def\KK{{\Bbb K}}

\def\PP{{\Bbb P}}
\def\QQ{{\Bbb Q}}

\def\ZZ{{\Bbb Z}}
\def\add{{\sanse add}}
\def\cov{{\sanse cov}}

\def\Cof#1{{{\sanse cof}$({\cal #1})$}}
\def\Unif#1{{{\sanse unif}$({\cal #1})$}}
\def\Cov#1{{{\sanse cov}$({\cal #1})$}}
\def\Add#1{{{\sanse add}$({\cal #1})$}}
\def\sm{{\smallskip}}
\def\ce#1{{\centerline{#1}}}
\def\no{{\noindent}}
\def\la{{\langle}}
\def\ra{{\rangle}}
\def\ha{{\hat{\;}}}
\def\sub{\subseteq}

\def\sem{\setminus}
\font\small=cmr8 scaled\magstep0
\font\smalli=cmti8 scaled\magstep0
\font\capit=cmcsc10 scaled\magstep0
\font\capitg=cmcsc10 scaled\magstep1

\font\dunhgg=cmdunh10 scaled\magstep2

\font\sanse=cmss10 scaled\magstep0
\font\bold=cmssdc10 scaled\magstep1
\font\bolds=cmssdc10 scaled\magstep0

\overfullrule=0pt
\openup1.5\jot

\ce{}
\Smallskip
\ce{\dunhgg Evasion and prediction II}
\footnote{}{{\openup-6pt {\small {\smalli
1991 Mathematics subject classification.}
03E05 03E35 20K20 20K25 \par
{\smalli Key words and phrases.} Specker
phenomenon, Cardinal invariants of the
continuum, Forcing, Evading and predicting
\endgraf}}}
\Bigskip
\no {\capit J\"org Brendle\footnote{$^1$}{{\small 
Supported by DFG--grant Nr. Br
1420/1--1.}}}
\sm
\item{} Mathematisches Institut der Universit\"at, Auf der Morgenstelle
10, 72076 T\"ubingen, Germany; email: {\sanse 
jobr@michelangelo.mathematik.uni--tuebingen.de}
\bigskip
\no {\capit Saharon Shelah\footnote{$^2$}{{\openup-7pt {\small 
Supported by the Edmund Landau
Center for research in Mathematical Analysis (sponsored by
the MINERVA--foundation (Germany)).\endgraf}}}$^,$\footnote{$^3$}{{\small 
publication
number 540 }}}
\sm
\item{} Department of Mathematics, Hebrew University, Givat Ram,
91904 Jerusalem, Israel
\par
\itemitem{} and
\par
\item{} Department of Mathematics, Rutgers University, New Brunswick,
NJ 08854, USA
\Bigskip

\ce{\capitg Abstract}
\bigskip
\no A subgroup $G \leq \ZZ^\omega$ exhibits the Specker phenomenon
if every homomorphism $G \to \ZZ$ maps almost all unit vectors
to $0$. We give several combinatorial characterizations of the
cardinal $\se$, the size of the smallest $G \leq \ZZ^\omega$
exhibiting the Specker phenomenon. We also prove the consistency 
of $\bb < \ee$, where $\bb$ is the unbounding number and $\ee$
the evasion number. Our results answer several questions
addressed by Blass.
\vfill\eject

{\bold Introduction}
\Smallskip
Specker [Sp] proved that given a homomorphism $h$ from $\ZZ^\omega
$ to the infinite cyclic group $\ZZ$, where $\ZZ^\omega$ denotes
the direct product of countably many copies of $\ZZ$,
we have $h(e_n) = 0$ for all but finitely many unit
vectors $e_n \in\ZZ^\omega$ (in other words, the $n$--th
component of $e_n$ is $1$, and its other components are $0$).
Blass [Bl] studied the {\it Specker--Eda number} $\se$,
the size of the smallest subgroup $G\leq\ZZ^\omega$ containing
all unit vectors which
still has the property that every homomorphism $h : G \to \ZZ$
annihilates almost all unit vectors. We will give various
(mostly less algebraic) characterizations of $\se$ (some
of which already play a prominent role in Blass' work);
we will also study some related cardinal invariants
of the continuum.
\par
To be more explicit, let $\leq^*$ denote the {\it eventual
domination order} on the Baire space $\omom$; i.e. $f \leq^*
g$ iff $f(n) \leq g(n)$ {\it for all but finitely many $n$}.
We shall usually abbreviate the statement in italics by
$\forall^\infty n$; similarly we will write $\exists^\infty n$
for {\it there are infinitely many $n$}. The {\it
unbounding number} $\bb$ is the smallest
size of a $\leq^*$--unbounded family ${\cal F}$ of functions
in $\omom$ (i.e., given any $g\in\omom$, there is $f\in{\cal F}$
with $\exists^\infty n \; (f(n) > g(n))$).
Given a $\sigma$--ideal ${\cal I}$ on $\omom$, the {\it additivity}
\add$({\cal I})$ is the least cardinality of a family
${\cal F}$ of members of ${\cal I}$ whose union is not in ${\cal I}$.
We shall use this cardinal only in the cases ${\cal I} = {\cal M}$,
the ideal of meager sets, and ${\cal I} = {\cal L}$, the ideal of
Lebesgue null sets. --- While the preceding invariants have been
studied by a number of people in the last two decades, the following
concept was introduced only recently by Blass [Bl].
Given an at most countable set $S$, an {\it $S$--valued predictor}
is a pair $\pi = (D_\pi , \la \pi_n ; \; n\in D_\pi \ra )$ where
$D_\pi \sub\omega$ is infinite and for each $n\in D_\pi$, $\pi_n $ is a function
from $S^n$ to $S$. $\pi$ {\it predicts } $f \in S^\omega$ iff for all but
finitely many $n \in D_\pi$, we have $f(n) = \pi_n (f \re n)$;
otherwise $f$ {\it evades} $\pi$. The {\it evasion number} $\ee$
is the smallest size of a family ${\cal F}$ of functions in $\omom$
such that no $\omega$--valued predictor predicts all $f\in{\cal F}$.
A $\ZZ$--valued predictor is {\it linear} iff all $\pi_n
: \ZZ^n \to \QQ$ are $\QQ$--linear maps. The corresponding
{\it linear evasion number} shall be denoted by $\ee_\ell$ (i.e.,
$\ee_\ell = \min \{ | {\cal F}| ; \; {\cal F} \sub \ZZ^\omega
$ and no linear $\ZZ$--valued predictor predicts all $f\in {\cal F}
\}$). (Blass' definition of linear evading [Bl, section 4] is slightly
different; however, it gives rise to the same cardinal; we use
the present definition because we shall work with functions
in $\ZZ^\omega$ in 2.2.)
\par
These notions enable us to phrase our main results.
\sm
{\bolds Theorem A.} {\it It is consistent with $ZFC$ to assume
$\bb < \ee$.}
\sm
{\bolds Theorem B.} {\it $\se = \ee_\ell = \min \{ \ee , \bb \}$.}
\sm
\no They will be proved in sections 1 and 2 of our work. Section 2 also
contains a further purely combinatorial characterization of the cardinal
$\se$ (subsections 2.4 and 2.5). To put our results into a somewhat
larger context, we point out the following consequences which
involve some earlier results, due mostly to Blass [Bl].
\sm
{\bolds Corollary.} {\it (a) \add$({\cal L}) \leq \se \leq$ \add$({\cal M})
\leq \bb$; \par
(b) any of the inequalities in (a) can be consistently strict; \par
(c) it is consistent with $ZFC$ to assume $\ee_\ell < \ee$.}
\sm
\no Theorems A and B together with the Corollary give a complete
solution to Questions (1) through (3) in [Bl, section 5].
Note in particular that the cardinals (2) through (5) in
Corollary 8 in [Bl, section 3] are indeed equal.
\sm
{\capit Proof of Corollary.} (a) This follows from Theorem B and
Blass' results [Bl, Theorems 12 and 13]. The well--known inequality
\add$({\cal M}) \leq\bb$ is due to Miller [Mi].
\par
(b) The consistency of \add$({\cal M}) < \bb$ is well--known
(it holds e.g. in the Mathias or Laver real models); for the
consistency of \add$({\cal L}) < \se$ see [Bl] (in particular
[Bl, Theorem 9]); the consistency of $\se < $ \add$({\cal M})$
follows from Theorem B and [Br, Theorem A].
\par
(c) This is immediate from Theorems A and B. $\qed$
\par\sm
A set of reals predicted by a single predictor is small in various senses;
it belongs, in particular, both to ${\cal M}$ and ${\cal L}$.
This motivates us to introduce the $\sigma$--ideal
${\cal J}$ on $\omom$ generated by such sets of reals
(see [Br, section 4] for more on this). Clearly, the uniformity
of ${\cal J}$ (i.e., the size of the smallest set of reals
not in ${\cal J}$) is closely related to the evasion number.
In fact, $\ee \leq \ee (\omega)$ where $\ee (\omega)$ denotes the
former cardinal. We shall show in section 3 that these two cardinals are
equal under some additional assumption, thus giving a partial 
answer to [Br, section 6, question (4)].
\par
The results of this work are due to the second author. It was
the first author's task to work them out and to write
up the paper.
\bigskip
{\bolds Notational remarks.} A p.o. $\PP$ is {\it $\sigma$--centered} iff
there are $\PP_n \sub \PP$ ($n\in\omega$) so that $\PP = \bigcup_n
\PP_n$ and given $n\in\omega$, $F \sub \PP_n$ finite, there 
is $p\in\PP$ extending all $q\in F$. $\PP$--names are denoted
by symbols like $\dot h , \dot\pi , \dot D $, ... $|$ stands
for {\it divides}; $\not |$ means {\it does not divide}.
\bigskip
{\bolds Acknowledgment.} The first author thanks the Wroc{\l}aw
logic group for several discussions about the material
presented here.
\Bigskip

{\bold $\S$ 1. Proof of Theorem A}
\Smallskip
1.1. We shall use a finite support iteration of $ccc$ p.o.'s
of length $\kappa$ (where $\kappa \geq \omega_2$ is a regular
cardinal) over a model $V$ for $CH$ 
to prove the consistency of $\ee > \bb$.
In fact, in the resulting model, $\bb = \omega_1$
and $\ee = \kappa$. We start with defining the p.o. $\PP$
we want to iterate. Notice that it is quite
similar to the one used in [Br, 4.3.] for predicting
below a given function.
\sm
\ce{$ \la d, \pi , F \ra \in \PP \Longleftrightarrow
d \in \twolom$ is a finite partial function,\hskip 5truecm}
\par
\ce{\hskip 3.5truecm $\pi = \la \pi_n ; \; n\in d^{-1} (\{ 1 \} )
\ra$ and $ \pi_n : \omega^n \to \omega$ is a finite partial function,}
\par
\ce{\hskip 3truecm $F \sub \omom$ is finite and $(f \neq g \in F 
\longrightarrow \max \{ n ; \; f\re n = g \re n \} < | d|)$.}
\sm
\no The order is given by:
\sm
\ce{$\la d' , \pi ' , F ' \ra \leq \la d , \pi , F \ra
\Longleftrightarrow d' \supseteq d , \pi ' \supseteq \pi ,
F ' \supseteq F $ and \hskip 5truecm}
\par
\ce{ \hskip 4truecm $(f \in F , n \in (d ' )^{-1} ( \{ 1 \} )
\sem d^{-1} ( \{ 1 \} ) \longrightarrow \pi_n ' (f \re n )
= f(n))$} 
\par
\ce{ \hskip 1truecm (in particular $\pi_n ' (f \re n)$ is defined).}
\sm
\no Notice that we use the convention that stronger
conditions are smaller in the p.o. --- The first two
coordinates of a condition are intended as a finite
approximation to a generic predictor; the third coordinate
then guarantees that functions are predicted
from some point on. Thus it is straightforward that
$\PP$ adjoins a predictor which predicts all ground--model
functions. Hence iterating $\PP$ increases $\ee$.
\par
Furthermore $\PP$ is $\sigma$--centered (and thus in particular
$ccc$). To see this simply notice that conditions with
the same initial segment in the first two coordinates
are compatible.
\par
So it remains to show that $\bb = \omega_1$ after iterating
$\PP$. For this it suffices to show the following:
\sm
\item{$(*)$} whenever $G \in W$ is an unbounded family
of functions from $\omega$ to $\omega$, and $\PP \in W$ is the p.o.
defined above, then
\par
\ce{$\forces_\PP `` G$ is unbounded".}
\sm
\no Using $(*)$ we can show that $\omom \cap V$ is still unbounded
in the final model: $(*)$ guarantees that it stays unbounded
in successor steps of the iteration; and one of the
usual preservation results for finite support iterations
(see, e.g., [JS, Theorem 2.2]) shows that it does so in limit steps of
the iteration as well. Now, $V \models CH$; hence $\omom \cap
V$ is an unbounded family of size $\omega_1$ in the final model.
\par
To start with the proof of $(*)$, let $ \dot h$ be a $\PP$--name
for a function in $\omom$.
For each $d \in \twolom , \pi = \la \pi_n ; \; n
\in d^{-1} (\{ 1 \} ) \ra$ an initial segment of a predictor
(as in the definition of $\PP$), $k \in \omega$ and 
$\bar f^* = \la f_\ell^* \in \omega^{|d|} ; \; \ell
< k \ra$  we define $h = h_{d, \pi , \bar f^*} \in (\omega + 1)^\omega$
by
\sm
\ce{$h(n) : = \min \{ m\leq\omega ; \; $ for no $p \in \PP$
with $p = \la d, \pi , F \ra$, $F = \{ f_\ell ; \;\ell < k \}
$, $f_\ell \re |d| = f_\ell^*$,}
\par
\ce{do we have $p \forces_\PP `` \dot h (n) > m " \}$.}
\bigskip
{\bolds 1.2. Main Claim.} {\it $h \in \omom$.}
\bigskip
{\capit 1.3. Proof of $(*)$ from the Main Claim.}
Let $h^* \in\omom$ such that for all $d , \pi , \bar f^*$
as above we have $h_{d , \pi , \bar f^*} \leq^* h^*$.
As $G$ is unbounded we can find $f\in G$ such that 
there are infinitely many $n$ with $f(n) > h^* (n)$.
We claim that $\forces_\PP ``\exists^\infty n (f(n) > \dot h(n) )"$.
This will show $(*)$.
\par
Assume $m\in\omega$ and $p\in\PP$ are such that
\sm
\ce{$p \forces_\PP `` \forall n \geq m \; (f(n) \leq \dot h (n))"$.}
\sm
\no Find $d, \pi , \bar f^*$ such that $p = \la d , \pi , F \ra$ where
$F = \{ f_\ell ; \;\ell < k \}$ and $f_\ell \re |d|
= f^*_\ell$. Find $n \geq m$ such that $f(n) > h^* (n)$ and
$h^* (n) \geq h_{d , \pi , \bar f^*} (n)$. Then
\sm
\ce{$p \forces_\PP `` h_{d , \pi , \bar f^*} (n)
< f(n) \leq \dot h (n) "$,}
\sm
\no contradicting the definition of $h_{d , \pi , \bar f^*}$.
$\qed$
\bigskip
{\capit 1.4. Proof of the Main Claim (1.2.).}
Let $d , \pi , k , \bar f^* = \la f_\ell^* ; \; \ell < k \ra
$ as above
and $n \in \omega$ be fixed.
Now assume that we have $p_i = \la d , \pi , \{ f_\ell^i
; \;\ell < k \} \ra$ with $f^i_\ell \re |d| = f^*_\ell$ and
\sm 
\ce{$p_i \forces_\PP `` \dot h (n) > i "$.}
\sm
\no We shall reach a contradiction.
As we can replace $\la p_i ; \; i\in\omega \ra$ by a 
subsequence, if necessary, we may assume that for all $\ell < k$:
\sm
\itemitem{{\sanse either}} $(a)_\ell$ for some $g_\ell \in
\omega^\omega$ for all $i$ $(f_\ell^i \re i = g_\ell \re i )$
\par
\itemitem{{\sanse or}} $(b)_\ell$ for some $i_\ell \in\omega$ and
$\hat g_\ell \in \omega^{i_\ell}$ $( f_\ell^i \re i_\ell = \hat g_\ell \;
\land\; f_\ell^i (i_\ell) > i )$.
\par
\sm
\no Notice that $i_\ell \geq |d|$ in the latter case. --- Let $d^* :=
d\cup 0_{[|d| , \max (i_\ell ; \; (b)_\ell \;{\rm holds}) + 1)}$;
i.e. the function $d^*$ takes value $0$ between $|d|$ and the maximum
of the $i_\ell$. Put $F^* := \{ g_\ell ; \; (a)_\ell$ holds$\}$.
Then clearly $p^* = \la d^* , \pi , F^* \ra \in \PP$. Now choose
$\ell^*$ and $q \leq p^*$ such that 
\sm
\ce{$q \forces_\PP `` \dot h (n)
= \ell^* "$.}
\sm
\no We shall find $i > \ell^*$ so that $q$ and $p_i$
are compatible; this is a contradiction because $q$ and $p_i$ force 
contradictory statements.
\par
Assume $q = \la d^q , \pi^q , F^q \ra$. Choose $i \geq \ell^*$
large enough such that:
\sm
\item{(A)} $i \geq |d^q|$;
\par
\item{(B)} $i \geq \max \{ \max \{ \sigma (j) ; \; \sigma\in dom
(\pi^q_m) \;\land\; j\in m \} ; \; m \in (d^q)^{-1} (\{
1\} ) \}$.
\par\sm
\no Notice that (A) implies that $f^i_\ell \re |d^q | = g_\ell \re
|d^q|$ whenever $(a)_\ell$ holds, while $f^i_\ell (i_\ell ) >
\max \{ \max$ $\{\sigma (j); \; \sigma\in dom (\pi^q_m ) 
\;\land\; j\in m\} ; \; m \in (d^q)^{-1} (\{ 1 \} ) \}$
by (B) in case $(b)_\ell$ holds. 
For such $i$ let $q^i = \la d^i , \pi^i , F^i \ra$ where
\sm
\item{---} $d^i = d^q \cup 0_{[ |d^q| , a) }$, where $a$ is large enough
such that all functions in $F^i$ disagree before $a$;
\par
\item{---} $\pi^i \supseteq \pi^q$ such that for all $m \in (d^q)^{-1}
(\{ 1\}) \sem d^{-1} (\{ 1 \} )$ 
and all $f^i_\ell$ so that $(b)_\ell$ holds, we have
\par
\ce{$f^i_\ell (m) = \pi^i_m (f^i_\ell \re m ). \;\;\;\;\; (\star)$}
\par
\item{} (This can be done because, by (B), $\pi_m^q$ was not defined yet 
for sequences of the form $f^i_\ell \re m$.)
\par
\item{---} $F^i = F^q \cup \{ f^i_\ell ; \; \ell < k \}$.
\par\sm
\no Now we clearly have $q^i \in \PP$ and $q^i \leq q$. So we are left
with checking $q^i \leq p_i$. The inclusion relations are all satisfied.
Hence it suffices to see that for all $\ell < k$ and $m \in (d^i)^{-1}
(\{ 1 \} ) \sem d^{-1} (\{ 1 \})$, we have
\sm
\ce{$ f^i_\ell (m) = \pi^i_m (f_\ell^i \re m ) . \;\;\;\;\; (+)$}
\sm
\no In case $(b)_\ell$ holds this is true by $(\star)$. In case
$(a)_\ell$ holds we have $f^i_\ell \re (m+1) = g_\ell \re (m+1)$
for all such $m$. As $q \leq p^*$ we have $\pi^i_m (g_\ell \re m) = \pi_m^q
(g_\ell \re m) = g_\ell (m)$ for such $m$, and $(+)$ holds again.
This completes the proof of the Main Claim. $\qed$
\Bigskip\vfill\eject

{\bold $\S$ 2. Proof of Theorem B}
\Smallskip
{\bolds 2.1. Theorem.} {\it $\se \leq \ee$.}
\sm
{\capit Proof.} Let ${\cal F} \sub \omom$, $| {\cal F} |
< \se$. By Blass' result $\se \leq \bb$ [Bl, Theorem 2], there is $g\in\omom$
such that for all $f\in{\cal F} \; \forall^\infty n \; (f(n) < g(n))$.
Without loss $g$ is strictly increasing. We let $\la p_n ; \;
n\in\omega\ra$ be a sequence of distinct primes such that $p_n >>
g(n)$ and $p_n >> \prod_{\ell < n} p_\ell$. For $f\in{\cal F}$,
let $a_f \in\omom$ be defined by
\sm
\ce{$ a_f (n) := f(n) \cdot \prod_{\ell \leq n} p_\ell$.}
\sm
\no Let $G \leq \ZZ^\omega$ be the pure closure of the subgroup
generated by the unit vectors $e_n$,
$n\in\omega$, and the $a_f$, $f\in {\cal F}$. Clearly $| G |
< \se$. Hence there is $h : G \longrightarrow \ZZ$ a homomorphism
such that $W: = \{ n ; \; h(e_n) \neq 0 \}$ is infinite.
\par
Let us define
\sm
\ce{$ W^* := \{ n \in\omega ; \; \exists i > n \; (p_i | h(e_m)$
whenever $m\in\{n+1 , ..., i-1\}$ but $p_i \not{|} h(e_n) ) \}$.}
\sm
\no We claim that $W^*$ is an infinite subset of $W$. To see this, first note
that trivially $W^* \sub W$, by the clause $p_i \not{|} h(e_n)$.
Next, given $n_0 \in W$, find $i > n_0$ so that $p_i \not{|}
h(e_{n_0})$. Then clearly there is $n \geq n_0$ so that $n \in W$ and 
$p_i \not{|} h(e_n)$ and for all $m \in \{ n+1 , ..., i-1 \}$,
$p_i | h(e_m)$. Thus $n \in W^*$. This shows that $W^*$ is infinite.
\par
We introduce a predictor $\pi = (W^* , \la \pi_n ; \; n \in W^* \ra )$
as follows. Given $n \in W^*$ and $s \in \omega^n$ so that $\max rng (s)
< g(n-1)$, if there is $f \in {\cal F}$ with $s \sub f$ and $f(n) < g(n)$
and $|h(a_f)| < p_{n-1}$, then let $\pi_n (s) = f(n)$ for some $f$
with the above property. Otherwise $\pi_n (s)$ is arbitrary.
\par
We claim that $\pi$ predicts all $f\in{\cal F}$. This clearly finishes
the proof. Assume this were false, i.e. there is $f\in{\cal F}$ which 
evades $\pi$. Let $n \in W^*$ be large enough, such that
$\max rng (f\re n) < g(n-1)$, $f(n) < g(n)$, $|h(a_f)| < p_{n-1}$ and
$\pi_n (f\re n) \neq f(n)$. Then, by the definition of $\pi$,
there must be $f' \in {\cal F}$ with $f'\re n = f \re n$, $f' (n)
< g(n) $, $|h(a_{f'})|
< p_{n-1}$ and $\pi_n (f' \re n ) = f' (n) \neq f(n)$.
Now, for $k \in \{ f , f' \}$, we let
$$\eqalign{a^0_k &= (a_k (0) , ... , a_k (n-1) , 0 , ...)\cr
a^1_k &= (0, ... , 0 , a_k (n) , 0 , ... ) \cr
a^2_k &= (0, ... , 0 , a_k (n+1) , ..., a_k (i-1) ,0 , ...)\cr
a^3_k &= (0, ... , 0 , a_k (i) , a_k (i+1) , ...)\cr}$$
where $i$ witnesses that $n \in W^*$. So we have $a_k = a_k^0
+ a^1_k + a^2_k + a^3_k$. Thus
$$ h(a_{f'} - a_f) = h(a^0_{f'} - a^0_f ) + h(a^1_{f'} - a^1_f ) +
h(a^2_{f'} - a^2_f ) + h(a^3_{f'} - a^3_f ) . \;\;\;\;\; (\star)$$
Clearly $h(a^0_{f'} - a^0_f ) = h(0) = 0$. Next, $p_i \cdot
\prod_{\ell \leq n} p_\ell$ divides $h(a^3_{f'} - a^3_f)$ by definition
of the $a_k$; it also divides $h(a^2_{f'} - a^2_f)$ by definition
of the $a_k$ and because $p_i | h(e_m)$ for $m \in \{ n+1 , ... , i-1 \}$
as $i$ witnesses $n \in W^*$. Thus $(\star)$ yields the equation
$$ h(a_{f'} - a_f) = h(a_{f'}^1 - a^1_f ) \;\;\;\;\; {\rm in }\;
\;\;\;\; \ZZ / (p_i \cdot
\prod_{\ell \leq n} p_\ell ) \ZZ .\;\;\;\;\; (\star\star)$$
The right--hand side in $(\star\star)$ must be non--zero,
because $p_i \not{|} h(e_n)$ (as $i$ witnesses $n \in W^*$) and
$p_i \not{|} (a_{f'} (n) - a_f (n)) = \prod_{\ell \leq n}
p_\ell \cdot (f' (n) - f(n) )$ (as $f' (n) , f(n) < g(n) << p_n
<< p_i$). However, it certainly is divisible by $\prod_{\ell \leq n}
p_n$, whereas the left--hand side in $(\star\star)$ is not unless it is
zero (as $|h(a_f)| , |h(a_{f'})| < p_{n-1} << p_n$). This shows that the
equation $(\star\star)$ cannot hold, the final contradiction.
$\qed$
\sm
Note that this result improves [Br, Theorem 3.2].
\bigskip
{\bolds 2.2. Lemma.} {\it $\ee_\ell \geq \min \{ \ee , \bb \}$.}
\sm
{\capit Proof.} Let ${\cal F} \sub \ZZ^\omega$, $| {\cal F}|
< \min \{ \ee , \bb \}$. Find $g \in\omom$ strictly increasing
so that for all $f \in {\cal F}$, we have $|f| <^* g$, where
$| f| (n) = |f(n)|$.
We partition $\omega$ into intervals $I_n$, $n \in\omega$,
so that $\max (I_n) + 1 = \min (I_{n+1})$, as follows.
$I_0 = \{ 0 \}$. Assume $I_n$ is defined; choose $I_{n+1}$
so that $| I_{n+1} | > [2 \cdot g( \max (I_n))]^{\sum_{i \leq n}
| I_i |}$. For $f \in {\cal F}$, define $\bar f$ by
$\bar f (n) : = f \re I_n$, and let $\bar{\cal F} = \{ \bar f ;
\; f\in {\cal F }\}$. Use $| \bar{\cal F} | < \ee$ to get a single
predictor $\bar \pi = (\bar D , \la \bar \pi_n ; \; 
n \in \bar D \ra )$ predicting all the $\bar f \in \bar{\cal F}$.
For $n \in \bar D$, let $\Gamma_n := rng (\bar \pi_n \re
(-g(\max(I_{n-1})) , g(\max(I_{n-1})))^{\bigcup_{i<n} I_i} )
\cap \ZZ^{I_n}$. So $|\Gamma_n | < |I_n |$; hence for some
$i_n \in I_n$, the vector $\bar x_{i_n} = \la t(i_n) ; \; 
t \in \Gamma_n \ra$ depends on the vectors $\{ \bar x_i =
\la t(i) ; \; t \in \Gamma_n \ra ; \; \min(I_n) \leq i < i_n
\}$. Say $\bar x_{i_n} = \sum_{\min(I_n) \leq i < i_n}
q_i^n \bar x_i$, where $q^n_i \in \QQ$. In particular,
for fixed $t \in \Gamma_n$, we have $t(i_n) =
\sum_{\min (I_n) \leq i < i_n} q_i^n t(i)$. This allows
us to define a linear predictor $\pi = ( D , \la \pi_n ;
\; n \in D \ra )$ with $D = \{ i_n ; \; n\in\omega\}$ and
$\pi_{i_n} (s) = \sum_{\min (I_n) \leq i < i_n }
q_i^n s(i)$. Note that if $n \in\omega$ is such that
$\max rng (|f| \re \cup_{i<n} I_i ) < g (\max (I_{n-1}))$ and
$\bar \pi_n (\bar f \re n) = \bar f (n)$, then $\pi_{i_n}
(f \re i_n) = f(i_n)$. Hence, as $\bar \pi$ predicts
all $\bar f \in \bar{\cal F}$, $\pi$ predicts all
$f \in {\cal F}$. $\qed$
\bigskip
2.3. Clearly, Theorem B follows from 2.1., 2.2. and Blass'
results $\ee_\ell \leq \se \leq \bb$ [Bl, Theorem 2, Corollary
8 and Theorem 10]. $\qed$
\bigskip
{\bolds 2.4. Definition.} Given $D \sub \omega$
infinite and $\bar a = \la a_n \in [\omega]^{\leq n} ; \; n \in D\ra$,
the {\it slalom} $S_D^{\bar a}$ is the set of all functions
$f$ in $\omom$ with $f(n) \in a_n$ for almost all
$n\in D$. 
\par
Using this notion we can give a combinatorial characterization
of the cardinal $\ee_\ell = \se$.
\bigskip
{\bolds 2.5. Lemma.} {\it $\min \{ \ee , \bb \} = \min
\{ | {\cal F} | ; \; {\cal F} \sub \omom$ and for all
$D \sub \omega$ and $\bar a = \la a_n \in [\omega]^{\leq n} ; \;
n \in D \ra$ there is $f \in {\cal F } \sem S_D^{\bar a} \}$.}
\sm
{\capit Note.} It is immediate that the cardinal on the
right--hand side is bigger than or equal to the additivity
of Lebesgue measure \add $({\cal L})$, by Bartoszy\'nski's
characterization of that cardinal ([Ba 1], [Ba 2]).
We also note that the original proof of \add$({\cal L}) \leq$
\add$({\cal M})$ [Ba 1] shows in fact that this
cardinal is $\leq$ \add$({\cal M})$ as well. This gives
an alternative proof of Blass' $\min \{ \ee , \bb \}
\leq$ \add$({\cal M})$ [Bl, Theorem 13].
\sm
{\capit Proof.} $" \geq "$. By Theorem B, it suffices to show
that $\ee_\ell$ is bigger than or equal to the cardinal
on the right--hand side. However, this is exactly like Blass' original proof
of \add $({\cal L}) \leq \ee_\ell$ [Bl, Theorem 12], and we therefore
leave details to the reader.
\par
$" \leq "$. This argument is very similar to the one in
Lemma 2.2. So we just stress the differences.
\par
Take ${\cal F }\sub\omom$, $|{\cal F}| < \min \{ \ee , \bb \}$.
Find $g$ strictly increasing and eventually dominating all functions from
${\cal F}$. As before, partition $\omega$ into intervals $I_n$,
$n\in\omega$; this time we require that $i_{n+1} : = g (\max (I_n))^{\sum_{
i\leq n} |I_i|} \in I_{n+1}$. $\bar f $, $\bar{\cal F}$ and $\bar\pi$,
$\bar D$ are defined as before. \par
We put $D:= \{ i_n ; \; n \in \bar D \}$ and $a_{i_n} = \{
\bar\pi_n (s) (i_n); \; s \in g (\max (I_{n-1}))^{\bigcup_{i<n} I_i}
\} \in [\omega]^{\leq i_n}$, and leave it to the reader to check
that ${\cal F} \sub S_D^{\bar a}$. $\qed$
\bigskip
2.6. The notion of linear predicting can be generalized as 
follows (see [Br, section 4] for details). Let $\KK$ be an
at most countable field. A $\KK$--valued predictor $\pi
= (D_\pi , \la \pi_n ; \; n\in D_\pi\ra )$ is {\it linear }
iff all $\pi_n : \KK^n \to \KK$ are linear. $\ee_\KK$
is the corresponding {\it linear evasion number}. We 
easily see $\ee_\QQ = \ee_\ell$. Rewriting the proof
of 2.2. in this more general context gives $\ee_\KK
\geq \min \{ \ee, \bb \}$ for arbitrary $\KK$ and
$\ee_\KK \geq \ee$ in case $\KK$ is finite. As $\ee_\KK
\leq \bb$ for infinite $\KK$ [Br, 5.4.], we get $\ee_\KK
=\min \{ \ee , \bb \}$ for such fields --- in particular
all $\ee_\KK$ for $\KK$ a countable field are equal. 
We do not know whether this is true for finite $\KK$.
Note that $\ee_\KK > \ee , \bb$ is consistent for
such fields [Br, section 4].
\Bigskip

{\bold $\S$ 3. Some results on evasion ideals}
\Smallskip
{\bolds 3.1. Definition.} We say a predictor $\pi = (D , \la \pi_n ;
\; n \in D \ra )$ {\it predicts} a function $f\in\omom$ {\it everywhere}
if $\pi_n (f \re n) = f(n)$ holds for all $n \in D$. We put $\ee
(\omega) := \min \{  |{\cal F}| ;\; {\cal F} \sub\omom \;\land\;
$ for all countable families of predictors $\Pi$ there is $f\in {\cal F}$
evading all $\pi\in\Pi\}$, the {\it uniformity of the evasion
ideal ${\cal J}$}. --- As usual, \cov$({\cal M})$ denotes the
{\it covering number} of the ideal ${\cal M}$, i.e. the smallest size
of a family ${\cal F} \sub {\cal M}$ so that
$\bigcup {\cal F} = \omom$.
\bigskip
{\bolds 3.2. Observation.} {\it Assume $\la D^n ; \; n \in \omega \ra$
is a decreasing sequence of infinite subsets of $\omega$, and $\la \pi^n
= (D^n , \la \pi^n_k ; \; k \in D^n \ra ) ; \; n\in\omega \ra$ is a
sequence of predictors. Then there are a set $D \sub \omega$, almost
included in all $D^n$, and a predictor $\pi = (D , \la \pi_k ;
\; k \in D \ra )$ predicting all functions which are predicted
by one of the $\pi^n$.}
\sm
{\capit Proof.} We can assume that each function which is predicted
by some $\pi^n$ is predicted everywhere by some $\pi^m$ ---
otherwise go over to sequences $\la E^n ; \; n \in\omega\ra$ and
$\la \bar\pi^n = (E^n , \la \bar\pi_k^n ; \; k \in E^n \ra ) ; \;
n\in\omega\ra$ such that (i) for all $n\in\omega$ there is $m\in\omega$
so that $E^m \sub D^n$ and $\bar \pi_k^m = \pi_k^n$ for $k
\in E^m$ and (ii) for all $n,m\in\omega$ there is $\ell\in\omega$
so that $E^\ell \sub E^n \sem m$ and $\bar\pi_k^\ell = \bar\pi_k^n$
for $k \in E^\ell$.
\par
Choose $d^n \in D^n$ minimal with $d^n > d^{n-1}$, and
put $D = \{ d^n ; \; n\in\omega\}$. Fix $n\in\omega$ and $s\in \omega^{d^n}$.
To define $\pi_{d^n} (s)$, choose, if possible, $i \leq n$ minimal
so that for all $k \in D^i \cap d^n$, we have $\pi^i_k (s \re k) = s(k)$,
and let $\pi_{d^n} (s) = \pi^i_{d^n} (s)$. If this is impossible,
let $\pi_{d^n} (s)$ be arbitrary.
\par
To see that this works, take $f\in\omom$ and $i \in \omega$ minimal
so that $\pi^i$ predicts $f$ everywhere. As the set of functions predicted
everywhere by a single predictor is closed, there are $n \geq i$ and 
$s \in \omega^{d^n}$ so that $s \sub f$ and $s$ is not predicted
everywhere by any of the $\pi^j$ where $j < i$. Then $\pi_{d^m} (f \re
d^m) = \pi^i_{d^m} (f \re d^m)$ for all $m \geq n$, as required.
$\qed$
\bigskip
{\bolds 3.3. Theorem.} {\it $\ee \geq \min\{ \ee (\omega) , $ \cov$({\cal M}) 
\}$; thus either $\ee <$ \cov$({\cal M})$ or $\ee (\omega) \leq$ 
\cov $({\cal M})$ imply
$\ee  = \ee (\omega)$.}
\sm
{\capit Remark.} The statement is very similar to a recent
result of Kamburelis who proved $\ss\geq \min \{ \ss (\omega) , $ \cov$({\cal 
M}) \}$, where $\ss$ is the splitting number and $\ss (\omega)$
the $\aleph_0$--splitting number.
\sm
{\capit Proof.} The second statement easily follows from the first.
To prove the latter,
let ${\cal F} \sub \omom$, $| {\cal F}| < \min \{\ee (\omega) , $ \cov$ 
({\cal M}) \}$. We shall show $| {\cal F} | < \ee$.
For $\sigma \in\omlom \sem \{\la\ra\}$, we construct recursively sets
$D^\sigma \sub\omega$ and predictors $\pi^\sigma = (D^\sigma ,
\la \pi^\sigma_n ; \; n \in D^\sigma \ra )$ such that:
\sm
\item{(i)} $D^{\sigma \re i} \supseteq D^\sigma$ for $i \in | \sigma |$;
\par
\item{(ii)} for all $f\in{\cal F}$ and all $\sigma \in\omlom$ there is
$i\in\omega$ so that $f$ is predicted by $\pi^{\sigma\ha\la i \ra}
$.
\par\sm
\no First construct $\pi^{\la i\ra} = (D^{\la i \ra} , \la \pi_n^{\la i \ra}
; \; n \in D^{\la i \ra} \ra )$ satisfying (ii) by applying $|
{\cal F} | < \ee (\omega)$. 
\par
To do the recursion, assume $\pi^\sigma = (D^\sigma , \la
\pi^\sigma_n ; \; n \in D^\sigma \ra )$ is constructed for some fixed
$\sigma \in\omlom$. Given $f \in\omom$, define $f^\sigma$ by:
\sm
\ce{$ f^\sigma (i) : = f  (k_i^\sigma )$,}
\sm
\no where $\{ k^\sigma_i ; \; i\in\omega \}$ is the increasing
enumeration of the set $D^\sigma$. Let ${\cal F}^\sigma = \{ f^\sigma ;
\; f \in {\cal F} \}$. Again we get $\omega$ many predictors $\bar
\pi^{\sigma \ha \la i \ra} = (\bar D^{\sigma \ha\la i\ra}
, \la \bar \pi_n^{\sigma\ha\la i\ra} ;\; n\in \bar D^{\sigma\ha\la i \ra} \ra
)$, $i\in\omega$, so that every $f^\sigma \in {\cal F}^\sigma$ is predicted
by some $\bar \pi^{\sigma\ha\la i\ra}$. Let
$D^{\sigma \ha\la i\ra} = \{ k^\sigma_j ; \; j \in \bar 
D^{\sigma \ha\la i\ra} \}$.
Fix $j \in \bar D^{\sigma \ha \la i\ra}$ and $s \in \omega^{k_j^\sigma}$.
Let $\bar s \in \omega^j$ be defined by $\bar s (\ell) = s (k_\ell^\sigma)$.
Put $\pi^{\sigma\ha\la i\ra}_{k_j^\sigma} (s) := \bar \pi_j^{\sigma
\ha\la i\ra} (\bar s)$. Now it is easy to see that $\pi^{\sigma \ha\la i \ra}
$ predicts $f$ whenever $\bar \pi^{\sigma\ha \la i\ra}$
predicts $f^\sigma$. Thus (i) and (ii) hold.
This completes the recursive construction.
\par
Given $f \in\omom$, let $T_f = \{ \sigma \in\omlom ; \;$ for all $i
\leq | \sigma | \; (\pi^{\sigma\re i}$ does not predict $f$
everywhere)$\}$. By the above construction, the sets $[T_f]$ are
nowhere dense for $f \in {\cal F}$. As $| {\cal F} | <$ \cov$ ({\cal M})$,
there must be $g\in\omom \sem \bigcup_{f\in{\cal F}} [T_f ]$. Now use
the Observation (3.2.) to construct a new predictor from the
$\la \pi^{g\re n} ; \; n\in\omega \ra$ which will predict all
$f\in{\cal F}$. $\qed$
\bigskip
3.4. It is unclear whether $\ee = \ee (\omega)$ can be proved in
$ZFC$. In view of Theorem 3.3 it seems reasonable to ask first
\sm
{\capit Question.} {\it Is $\ee >$ \cov$({\cal M})$ consistent?}
\sm
\no Of course, we may also consider the cardinal $\ee_\ell (\omega)$,
the smallest size of a family ${\cal F}$ of functions from $\omega$
to $\omega$ such that no countable family of linear predictors
predicts all $f\in{\cal F}$. However, it is now easy to see that
$\ee_\ell (\omega) = \ee_\ell$. This is so because $\ee_\ell (\omega)
\leq \min \{ \ee (\omega) ,\bb\} \leq \min \{ \ee , \bb \} \leq
\ee_\ell$. To see the first inequality, note that the argument
for $\ee_\ell \leq \bb$ gives $\ee_\ell (\omega) \leq \bb$
as well (see [Br, section 5.4] for a stronger result);
for the second inequality, $\min \{ \ee (\omega) , \bb \} \leq$
\cov$({\cal M})$ by rewriting Blass' $\min \{ \ee , \bb \}
\leq $ \cov$({\cal M})$ [Bl, Theorem 13] and thus $\min \{ \ee
(\omega) , \bb \} = \min \{ \ee (\omega) , $ \cov$({\cal M}) , \bb \}
\leq \min \{ \ee , \bb \}$ by Theorem 3.3; the third
inequality is Lemma 2.2.
\bigskip
{\bolds 3.5. Duality.} Most of the cardinal invariants of the continuum
come in pairs and results about them usually can be dualized (see
[Br, section 4.5] for details). In our situation, the dual 
cardinals are: the {\it dominating number} $\dd$ (dual to $\bb$),
the smallest size of a family ${\cal F} \sub \omom$ such that
given any $g\in\omom$ there is $f\in{\cal F}$ with $g \leq^* f$;
the
{\it (linear) covering number} \cov$({\cal J})$ (\cov$({\cal J}_\ell)$)
of the ideal ${\cal J}$ (${\cal J}_\ell$) (the first being dual to
both $\ee$ and $\ee (\omega)$, the second being dual to $\ee_\ell$),
the least cardinality of a family of (linear) predictors $\Pi$
such that every function $f\in\omom$ ($\ZZ^\omega$) is predicted by some
$\pi\in\Pi$. Then we get:
\sm
{\capit Theorem.} {\it (a) It is consistent with $ZFC$ to assume
$\dd >$ \cov$({\cal J})$.
\par
(b) \cov$({\cal J}_\ell) = \max \{$ \cov$({\cal J}) , \dd \} = 
\min\{ | {\cal S}| ; \; {\cal S}$ consists of slaloms $S_D^{\bar a}$
where $\bar a = \la a_n \in [\omega]^{\leq n} ; \; n\in D \ra$
and $D \sub\omega$ is infinite and $\forall f\in\omom \; \exists
S_D^{\bar a} \in {\cal S} \;
\forall^\infty n \in D \; (f(n) \in a_n ) \}$. }
\sm
{\capit Proof.} These dualizations are standard, and we therefore
refrain from giving detailed proofs. The model for (a) is gotten
by iterating the p.o. $\PP$ from $\S$ 1 $\omega_1$ times
with finite support over a model for $MA + \neg CH$. (b) is the
dual version of Theorem B and Lemma 2.5. $\qed$
\sm
We close our work with a diagram showing the relations between
the cardinal invariants considered in this work (in particular,
the Specker--Eda number $\se$ and the evasion number $\ee$)
and some other cardinal invariants of the continuum (in
particular, those of Cicho\'n's diagram). We refer the reader
to [Bl], [Br] or [Fr] for the cardinals not defined here.
A similar diagram was drawn in [Br, section 4].
\bigskip
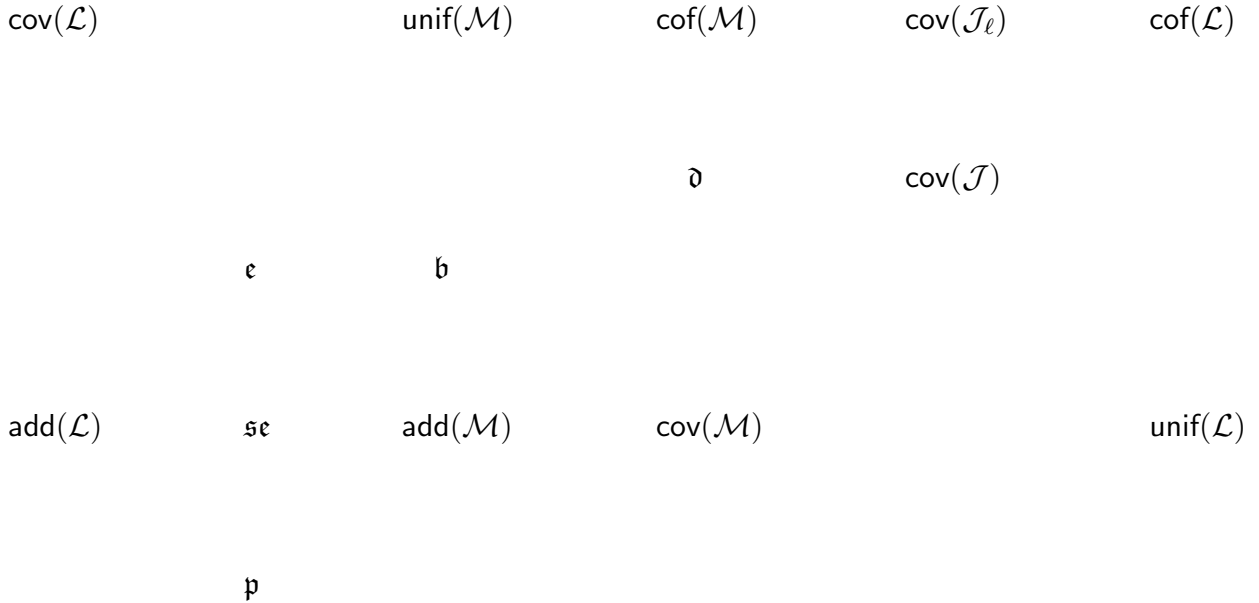
\halign{#  &\Hoskip # &\Hoskip  # &\Hoskip  # &\Hoskip  # &\Hoskip  # \cr
\Cov L && \Unif M & \Cof M & \cov$({\cal J}_\ell)$ & \Cof L \cr
\noalign{\Smallskip}
&  &  & \quad \dd & \Cov J &\cr
\noalign{\bigskip}
& \ee & \quad \bb &&&\cr
\noalign{\Smallskip}
\Add L & \se & \Add M & \Cov M && \Unif L \cr
\noalign{\Smallskip}
& \pp &&&& \cr}
\bigskip
\no In the diagram, cardinals increase as one moves up and to
the right. To enhance readability, we omitted the
relations $\ee \leq$ \Unif L, and its dual \Cov L $\leq
$\Cov J.
The dotted lines give the relations \Add M $=\min \{ \bb,$\Cov M $\}$,
$\se = \min \{ \ee , \bb \}$, and their dual versions.

\Bigskip

{\bold References}
\Smallskip
\itemitem{[Ba 1]} {\capit T. Bartoszy\'nski,} {\it Additivity of measure
implies additivity of category,} Transactions of the American Mathematical
Society, vol. 281 (1984), pp. 209-213. 
\smallskip
\itemitem{[Ba 2]} {\capit T. Bartoszy\'nski,} {\it Combinatorial
aspects of measure and category,} Fundamenta Mathematicae, vol. 127
(1987), pp. 225-239.
\smallskip
\itemitem{[Bl]} {\capit A. Blass,} {\it Cardinal characteristics
and the product of countably many infinite cyclic groups,}
Journal of Algebra (to appear).
\smallskip
\itemitem{[Br]} {\capit J. Brendle,} {\it Evasion and prediction ---
the Specker phenomenon and Gross spaces,} Forum Mathematicum
(to appear).
\sm
\itemitem{[Fr]} {\capit D. Fremlin,} {\it Cicho\'n's diagram,}
S\'eminaire Initiation \`a l'Analyse (G. Choquet,
M. Rogalski, J. Saint Raymond), Publications Math\'ematiques
de l'Universit\'e Pierre et Marie Curie, Paris, 1984,
pp. 5-01 - 5-13. 
\smallskip
\itemitem{[JS]} {\capit H. Judah and S. Shelah,} {\it The Kunen-Miller
chart (Lebesgue measure, the Baire property, Laver reals and
preservation theorems for forcing),} Journal of
Symbolic Logic, vol. 55 (1990), pp. 909-927.
\smallskip
\itemitem{[Mi]} {\capit A. Miller,} {\it Some properties of
measure and category,} Transactions of the American Mathematical
Society, vol. 266 (1981), pp. 93-114.
\smallskip
\itemitem{[Sp]} {\capit E. Specker,} {\it Additive Gruppen von
Folgen ganzer Zahlen,} Portugaliae Mathematica, vol. 9 (1950),
pp. 131-140.

\vfill\eject\end